\documentclass[12pt]{article}
\usepackage[ansinew]{inputenc}
\usepackage[pdftex]{graphicx}
\usepackage{amssymb,amsmath,amscd}
\usepackage{dsfont}
\usepackage{verbatim}
\usepackage{amsmath} 
\usepackage{graphicx}
\usepackage{pstricks}
\usepackage{enumitem}
\usepackage{ntheorem}
\usepackage{amssymb}
\usepackage{float}
\usepackage{tikz-cd}
\usepackage[left=2cm,top=2cm,right=2cm,bottom=2cm,nohead]{geometry}

\newtheorem{de}{Definition}[section] 
\newtheorem{nott}{Notation}

\newtheorem{theo}{Theorem}
  
\newtheorem{lem}{Lemma}[section]
 
\newtheorem{conj}{Conjecture}[section]

\setenumerate[1]{font=\bfseries,label=\Roman*.}
\setenumerate[2]{font=\itshape,label=\arabic*)}

\author{Yashar Memarian
}       
\title{A Note on the Geometry of Positively-Curved Riemannian Manifolds\\
       }
\date{}

\begin{document}
\maketitle
\tableofcontents

\begin{abstract}
In this paper I present a comparison theorem for the waist of Riemannian manifolds with positive sectional curvature. The main theorem of this paper gives a partial positive answer to a conjecture formulated by M.Gromov in \cite{grwst}. The content of this paper combines two aspects: classical volume comparison theorems of Riemannian geometry, and geometric measure theoretic ideas stemming from Almgren-Pitts Min-Max theory. 
\end{abstract}

\section{Introduction}
A basic way of studying or understanding the geometry and/or topology of a Riemannian manifold (or more generally a metric space) is by estimating metric invariants on such a space. But what is a metric invariant? For the class of Riemannian manifolds, we all know at least two basic metric invariants which are the volume and the diameter. Probably most of us are more familiar with the term \emph{topological invariant}. A topological invariant is an object $A$ (a number, a group, a ring, etc.) associated to any toplogical space $X$ such that any homeomorphism from $X$ to $X'$ induces an \emph{isomorphism} from $A$ to $A'$. In order to give a geometric flavour to the definition of a topological invariant, we could make the homeomorphism an isometry. Therefore, a metric invariant attached to a metric space $M$ is an object $A$ such that any \emph{isometry} between $M$ and $M'$ maps $A$ \emph{isomorphically} to $A'$. It is clear now why diameter and volume are considered \emph{metric invariants}. The problem with isometric maps is that they are very restrictive. Therefore it is best to look at a wider class of mappings. One way (and the way I shall utilise here) is to look at the class of $\lambda$-bi-Lipschitz maps which I will define now:
\begin{de}[$\lambda$-bi-Lipschitz maps]
A map $f:X\to Y$ from a metric space $(X,d)$ to the metric space $(Y,d')$ is said to be $\lambda$-bi-Lipschitz for $\lambda\geq 0$ if for every $x,y\in X$ we have
\begin{eqnarray*}
\frac{1}{\lambda}d(x,y)\leq d'(f(x),f(y))\leq \lambda d(x,y).
\end{eqnarray*}
\end{de}
Now we ask our metric invariants to behave \emph{as desired} under Lipschitz mappings. For example:
\begin{de}[Metric Invariants] \label{mi}
A \emph{metric invariant} is a map which associates to every metric space $X$ a non-negative number $Inv(X)\in\mathbb{R}_{+}$ such that any $\lambda$-bi-Lipschitz homeomorphism $f$ from $X$ to $X'$ map $Inv(X)$ to $Inv(X')$ such that
\begin{eqnarray*}
\frac{1}{F(\lambda)}Inv(X)\leq Inv(X')\leq F(\lambda)Inv(X),
\end{eqnarray*}
for a certain function $F:\mathbb{R}_{+}\to\mathbb{R}_{+}$.
\end{de}

\emph{Remarks}:
\begin{itemize}
\item One could define more complicated metric invariants. For example, if a class of (abelian) groups (or rings, fields, etc.) is such that each element of this class is endowed with a norm, then we can associate a normed (abelian) group to each metric space (defined to be a metric invariant) if under $\lambda$-bi-Lipschitz homemorphisms of spaces the induced maps will be $F(\lambda)$-bi-Lipschitz isomorphisms (considered with the norms of these groups).
\end{itemize}

Over the past fifty years, many metric invariants have been presented. In fact, defining metric invariants and estimating them on different geometric spaces composed a large part of Gromov's research career (for example see \cite{grbook}). Usually an exact estimation of a metric invariant for a given (even well-known) metric space is an extremely difficult task (for example, remember the filling volume of a circle). Hence, a metric invariant by itself is not of great use until one can compare it with other manifolds, called the model manifolds. Comparison geometry has a great history (at least in Riemannian Geometry) and more recently for more general metric spaces (for example the $cat(\kappa)$ spaces and particularily hyperbolic groups of geometric group theory). For more on comparison geometry the reader can consult the excellent books \cite{berger} and \cite{cheeger}. In this paper I am mainly interested in Riemannian manifolds with sectional curvature (denoted by $K$) everywhere $\geq \delta >0$. In positive curvature, there are lots of interesting open questions regarding the topology of such manifolds. With strictly positive curvature, I should point out that we do not have many examples at our disposition. An ambitious and interesting program of research in Riemannian Geometry of today is the construction of new examples of (strictly) positively-curved Riemannian manifolds. To do so, usually one begins by assuming that the manifold has a \emph{large} isometry group, then progressively reduces the \emph{size} of the isometry group.

In this paper, I will present a relatively new metric invariant which was defined by Gromov in \cite{grwst}, and prove a comparison type theorem for this invariant for positively curved Riemannian manifolds.

However before doing this, I find it useful to provide a recollection of some theorems concerning the comparison/estimation of several important metric invariants on Riemannian manifolds of positive (and in some case of non-negative) curvature.

The first theorem by F.Wilhelm and proved in \cite{will} involves the author comparing the filling radius of positively-curved Riemannian manifolds with the sphere of constant curvature. Recall the
\begin{de}[Filling Radius]
Let $M$ be a closed Riemannian manifold of dimension $n$ and $i:M\to X$ be an isometric embedding where $X$ is a metric space. Let $(M+\varepsilon)_{X}$ be the $\varepsilon$-neighborhood of $I(M)$ in $X$. The inclusion map $j: I(M)\to (M+\varepsilon)_{X}$ induces a map on the $n$-homology which is denoted by $i^{*}_{\varepsilon}:H_n(I(M))\to H_n((M+\varepsilon)_{X})$. The filling radius of $M$ relative to $I$ is the infimum of $\varepsilon$ such that $i^{*}_{\varepsilon}=0$. The filling radius of $M$ is defined to be relative to the embedding of $M$ into $L^{\infty}(M)$ where each point is mapped to the distance function to that point.
\end{de}
Intuitively, when our manifold is isometrically embedded into a metric space, begin enlarging your manifold until you trivialise the homology (equivalently until you bound something). This first enlargment value is then the filling radius relative to your embedding. The coefficient of the homology (depending on orientability) of the manifold is taken to be the set of integers for orientable manifolds, or integers modulo $2$ in case the manifold is non-orientable.
Then
\begin{theo}[Wilhelm]
Let $M$ be a closed Riemannian manifold of dimension $n$ having sectional curvature everywhere $\geq 1$ and $\mathbb{S}^n$ be the canonical Riemannian sphere (of constant curvature $1$. Then
\begin{eqnarray*}
FillRad(M)\leq FillRad (\mathbb{S}^n).
\end{eqnarray*}
\end{theo}
In fact, Wilhelm proves more. He shows that if the filling radius of $M$ is equal the filling radius of the sphere, then $M$ is isometric to the sphere, and if the filling radius is larger than $\pi/6$, then $M$ is a twisted sphere. He also proves the existence of a constant $C(n)$ (depending only on the dimension) that if 
\begin{eqnarray*}
FillRad(\mathbb{S}^n)-C(n)<FillRad(M)
\end{eqnarray*}
Then $M$ is diffeomorphic to $\mathbb{S}^n$.

Now we turn our attention to another (set of) invariants which are defined very similarily to the invariant waists defined in the next section. These invariants are called the $k$-widths, (see \cite{grwidth})
\begin{de}[k-Widths]
Let $X$ be a metric space. The $k$-dimensional width of $X$ denoted by $W_k(X)$ is the exact lower bound of $\delta>0$ for which there exists a $k$-dimensional space $T$ and a continuous map $f:X\to T$ for which all the fibers (inverse images) have diameter at most $\delta$.
\end{de}
The widths are defined through a min-max procedure. You look at the infimum of the set of \emph{diameters} obtained by projecting your space to $k$-dimensional (topological) space. If you study the widths of a subspace of a fixed metric space, you obtain the $k$-diameters. The idea of the definition of widths is to universalise the diameters (quite similar to what we said about filling radius). The difference is, that for a filling radius, Gromov shows in \cite{fillrad} that the best (or the least) filling radius is relative to $L^{\infty}$ but for widths such universal space is not known (to my knowledge).
 Be aware of the duality between the dimension and codimension, intuitively, the $k$-widths tell you how fat a space is (in the sense of diameter) in \emph{co}-dimension $k$ (since you look at pull-backs over a space of dimension $k$).
 
 The next theorem was conjectured by M.Gromov in \cite{grwidth} and proved in \cite{perel} by G.Perelman.
 \begin{theo}[Perelman]
 Let $M$ be a closed non-negatively curved Riemannian manifold. A constant $C(n)$ exists such that
 \begin{eqnarray*}
 C(n)^{-1}vol_n(M)\leq \displaystyle\Pi_{k=0}^{n-1}w_k(M)\leq C(n)vol_n(M).
 \end{eqnarray*}
 \end{theo}
 Remark that Perelman's theorem asserts the fact that non-negatively curved Riemannian manifolds \emph{metrically} look like $n$-dimensional rectangles where each side of the rectangle corresponds to the widths of the Manifold.
 Perelman's proof begins by proving a relaxed version of the assertion of the theorem, where one replaces the widths by a close (but easier to evaluate) invariant called the \emph{packing widths}. He proves the theorem using the collapse technique for packing widths and finishes the proof by geometrically relating the widths and packing widths. The sharp value of the constant is unknown.
 
Wilhelm and Perelman's theorems give satisfactory comparison/estimation theorems for the corresponding invariants. The purpose of this paper is to estimate another important metric invariant for positively-curved Riemannian manifolds.

\section{The Space of Cycles and the Invariant Waist}

The invariant waist is another important invention of M.Gromov in \cite{grwst}, and its history goes back to F.Almgren and his Morse Theory on the space of cycles (I will come to Almgren's Theory in detail later on this paper). Intuitively, the $k$-waists tell you how fat your space is (in the sense of volume this time) in codimension $k$. Because of all the technicalities due to the volume in lower (co)dimensions, with a bit more care I could say for the $k$-waists to take the definition of the widths and replace the diameters of the fibers by their Hausdorff measure and the same Min-Max processes to give you the waist. Of course, one needs to work in the class of metric-measure spaces rather than only metric spaces.
In order to define the invariant waist, I shall need to give the definition of the space of cycles, a space on which a variational procedure (min-max procedure) will be used.

\begin{de}[The space of cycles]
Let $X$ be a metric space and let $G$ be the group of coefficent. The space of $k$-cycles with coefficients in $G$, denoted by $C_k (X,G)$ is the set constituted by the singular Lipschitz $k$-cycles with coefficients in $G$. A singular Lipschitz simplex $\sigma:\Delta\to X$ possesses a volume. 
If $G$ is equipped with a norm, the mass of a $k$-cycle $T=\sum g_i \sigma_i$ is defined to be
\begin{eqnarray*}
M(T)=\sum \vert g_i \vert vol(\sigma_i),
\end{eqnarray*}
and the $\flat$ norm is defined by
\begin{eqnarray*}
\flat(T)=\inf\{M(S)+M(R)\,\vert\,T=S+\partial R\}.
\end{eqnarray*}
We equip $C_k (X,G)$ with the topology of the $\flat$ norm.
\end{de}

\emph{Remark}: Other topologies may be considered on $C_k(X,G)$. I will come back to this in the appropriate time.

Space of $k$-cycles being defined, we need to specify \emph{subsets} of this space which have significant geometric properties and will be important for our study:

\begin{de}[Parametrised Family of Cycles]
Let $T$

be an $m$-dimensional topological space. We say the space of $k$-cycles (of a metric space $X$) is parametrised by $T$ 
if there exists a continuous map $F$ from $T$ to $C_k(X,G)$.
\end{de}

Among the parametrised families of cycles defined on $X$, some particular ones are of interest which geometrically correspond to a \emph{sweep out} of the metric space $X$.
\begin{de}[Sweep Out]
Let $X$ be an $n$-dimensional metric space. $T$ an $(n-k)$-dimensional topological space. A family of $k$-cycles parametrised by $T$ and sweeping out the metric space $X$ corresponds to a continous map $F:T\to C_k(X,G)$ such that 
\begin{eqnarray*}
F^{*}a\neq 0,
\end{eqnarray*}
where $a\in H^{*}(C_k(X,G),G)$ is the \emph{fundamental} cohomology class of $C_k(X,G)$.
\end{de}

\emph{Remark}: Suppose $T$ is an $(n-k)$-dimensional simplicial complex (geometric measure theorisist prefers a cubical complex, both would do). The family $F$ corresponds to a sweep out of the metric space $X$, if (roughly) the following scenario happens:

The map $F$ sends every vertex of the simplicial complex $T$ to a $k$-cycle. For every two such vertices $c_1^0$ and $c_2^0$ which are the boundary of a $1$-simplex $c^1$, the image of this simplex by $F$, i.e. $F(c^1)$ is a $k+1$-dimensional chain such that:
\begin{eqnarray*}
\partial F(c^1)=F(c^0_1)-F(c^0_2).
\end{eqnarray*}
If the map $F$ preserves this pattern inductively on the dimension, then the family $F$ sweeps out the metric space by a family of $k$-cycles.

\emph{Remark}: This definition is different from the original definition of \emph{sweep out} proposed by Almgren and discussed in \cite{pitts}. Originally, sweep out is defined using the homotopy group of the space of cycles. We shall come back to this in Section $5$.

\begin{de}[k-Waists]
Let $X$ be a metric-measure space of dimension $n$. The $k$-waist of $X$ denoted by $wst_k(X)$ is the infimum of numbers $r\geq 0$ such that for every family of $k$-cycles (or relative cycles) parametrised by a $n-k$ dimensional $\mathbb{Z}_2$-topological manifold and generating the fundamental $\mathbb{Z}_2$-homology class of the space of cycles, the $k$-volume (Hausdorff measure of dimension $k$) of every cycle is at most equal to $r$.
\end{de}
One could avoid the topological complications of the space of cycles and define the waist relative to mappings in an easier way. In fact this second definition would be enough for the sake of this paper.
\begin{de}[k-Waists relative to spaces]
Let $X$ be a metric-measure space of dimension $n$. Let $T$ be a $(n-k)$-dimensional space. The $k$-waist of $X$ relative to $T$ , denoted by $wist_k(X\to T)$ is the exact number $r\geq 0$ such that for every continuous map $f:X\to T$, a $t\in T$ exists such that
\begin{eqnarray*}
Haus_k(f^{-1}(t))\geq r.
\end{eqnarray*}
\end{de}
Note that a variation of the above definition where instead of the $k$-volumes we consider the measure of the $\varepsilon$-neighborhood of the fibers could also be defined. To see this first set the following:
\begin{nott}
Let $X$ be a general metric-measure space and $Y\subset X$. Then
\begin{eqnarray*}
Y+\varepsilon =\{x \in X \vert \hspace{0.5mm} d(x,Y)\leq \varepsilon\},
\end{eqnarray*}
where $d(.,.)$ stands for the metric of $X$ and $d(x,Y)=inf_{y\in Y} d(x,y)$.
\end{nott}

I can now give another definition for the invariant waist.

\begin{de}[Waist of a mm-space, the variational viewpoint]
Let $X$ be a mm-space. Let $Z$ be a topological space and let $F_z$, $z \in Z$, be a family of cycles (subspaces) of $X$. Let $w(\varepsilon)$ be a function on $\mathbb{R}_+$. We say that the waist of $X$ relative to the family of cycles $F_z$ is at least equal to $w(\varepsilon)$, and we write
\begin{eqnarray*}
wst(X,F_z,\varepsilon) \geq w(\varepsilon),
\end{eqnarray*}
if there exists $z \in Z$ such that for every $\varepsilon>0$, 
\begin{eqnarray*}
\mu(F_z+\varepsilon) \geq w(\varepsilon).
\end{eqnarray*}
\end{de}

I recall Gromov's Conjecture stated in \cite{grwst}:
\begin{conj}[Gromov] \label{gro}
Let $X$ be an $n$-dimensional Riemannian manifold with sectional curvature $K\geq 1$. Let $f:X\to \mathbb{R}^k$ be a continous map. Then there is a $z\in\mathbb{R}^k$ such that for every $\varepsilon>0$ we have:
\begin{eqnarray*}
\frac{vol_n(f^{-1}(z)+\varepsilon)}{vol_n(X)}\geq \frac{vol_n(\mathbb{S}^{n-k}+\varepsilon)}{vol_n(\mathbb{S}^n)},
\end{eqnarray*}
where $\mathbb{S}^n$ is the canonical Riemannian $n$-sphere with sectional curvature equal to $1$.
\end{conj}
Recall that Conjecture \ref{gro} is proved for the case where $X$ itself is the canonical Riemannian sphere (see \cite{grwst} and \cite{memwst}). This conjecture hence presents a waist comparison for all Riemannian manifolds with $K\geq 1$.

The exposition of this paper will be as follows: at first, I prove a generalised Gromov Conjecture for the case $\varepsilon=0$, more precisely I will prove the following:

\begin{theo}\label{main}
Let $X$ be a closed Riemannian manifold of dimension $n$ (with possibly a non-empty quasi-regular convex boundary) such that we have $K\geq \lambda>0$ everywhere. Let $f:X\to \mathbb{R}^k$ be a smooth map. Then there exists a $z\in \mathbb{R}^k$ such that :
\begin{eqnarray*}
\frac{vol_{n-k}(f^{-1}(z))}{vol_n(X)}\geq \frac{vol_{n-k}(\mathbb{S}^{n-k}(\lambda))}{vol_n(\mathbb{S}^n(\lambda))}.
\end{eqnarray*}
where $vol_{n-k}$ stands for the Riemannian volume (or equivalently the Hausdorff measure) in dimension $(n-k)$ and $\mathbb{S}^n(\lambda)$ is the sphere of constant curvature equal to $\lambda$.
\end{theo}

The proof of this theorem relies on combining classical comparison type theorems in Riemannian Geometry and the Almgren-Pitts Min-Max theory. The next sections of this paper will concern this.

In the last section, I return to Conjecture \ref{gro} for $\varepsilon>0$, where possible ideas for proving this conjecture (at least for small enough $\varepsilon>0$) will be presented.


Remark that Conjecture \ref{gro} asserts that:
\begin{eqnarray*}
\frac{wst_k\{X\to\mathbb{R}^k,\varepsilon\}}{vol_n(X)}\geq \frac{wst_k\{\mathbb{S}^n(\lambda)\to\mathbb{R}^k,\varepsilon\}}{vol_n(\mathbb{S}^n(\lambda))}.
\end{eqnarray*}

I begin by reviewing some classical results in Riemannian geometry:

\section{Lower Bound on the Volume of Submanifolds with an Upper Bound on its Mean-Curvature}

This section concerns a few classical (volume) comparison-type results in Riemannian Geometry. One can find all the assertions of this section in advanced text books on Riemannian Geometry. However, to be thorough, I shall give the details of the computations here. I begin by setting up a few definitions and notations. All throughout this section (and this paper) I will be interested in the behaviour of the volume in a neighborhood of a submanifold of a (general) Riemannian manifold. The reasonable setting for studying such problems is the study of the normal exponential map with respect to a submanifold. In terms of coordinates, this leads to the study of Fermi coordinates with respect to a submanifold. Such coordinates are a generalisation of the normal coordinates, which is very common for studying the geometry of balls in Riemannian manifolds.

From now on, $M^n$ will be a closed $n$-dimensional manifold and $H^k$ a closed $k$-dimensional (smooth) submanifold. The normal bundle of $H$ is defined as follows:
\begin{eqnarray*}
N_H=\{(h,u)\vert h\in H, u\in TH_{h}^{\perp}\},
\end{eqnarray*}
where $TH_{h}^{\perp}$ denotes the orthogonal complement of the tangent space $TH_h$ in the tangent space $TM_h$. $N_H$ is a vector bundle over $H$ hence a differentiable manifold. The exponential map of the normal bundle $N_H$ is the map denoted by $exp_N$ defined by:
\begin{eqnarray*}
exp_N(h,u)=exp_h(u),
\end{eqnarray*}
where $(h,u)\in N_H$ and where $exp_h$ is the usual exponential map of $M$ at the point $h$. Hence we get a map:
\begin{eqnarray*}
exp_N:N_H\to M.
\end{eqnarray*}
It is intuitively clear (and not rigourously hard) to prove that $exp_N:N_H\to M$ maps a neighborhood of $H\subset N_H$ diffeomorphically onto a neighborhood of $H\subset M$. Define a \emph{tube} of radius $r\geq 0$ around $H$, denoted by $U_r(H)$ to be the set of all points $x\in M$ such that a geodesic $\gamma$ of length $\leq r$ from $x$ exists and meet $H$ orthogonally. Note that in general $U_r(H)\neq H+r$, where $H+r$ was defined previously and stands for the $r$-tubular neighborhood of $H$.

One should notice that $N_H$ inherits of a (canonical) Riemannian metric (hence a Riemannian volume). Indeed in order to define a metric on $N_H$, let $y=(p,v)\in N_H$ and let $w$ be a vector in the vector space $T_y(N_H)$. This means that $w$ is the tangent vector $\omega'(0)$ of a certain curve $\omega:t\to (p(t),v(t))$ with origin $y=\omega(0)$. Then defining
\begin{eqnarray*}
\vert w\vert^2=\vert\frac{\nabla p(0)}{dt}\vert^2+\vert\frac{\nabla v(0)}{dt}\vert^2
\end{eqnarray*}
defines a Riemannian metric on $N_H$ where $\nabla$ is the covariant derivative in $M$. Hence the volume element in $N_H$ can be represented in the form
\begin{eqnarray*}
dV_{N_H}(p,tv)=t^{n-k-1}dV_{H}(p)\,du(v)\,dt,
\end{eqnarray*}
where $\vert v\vert=1$ and $du$ is the volume element in the unit sphere $\mathbb{S}^{n-k-1}$ and $dV_{H}$ is the volume element in $H$.

The geometry of the submanifold $H$ is encoded in a bilinear form (tensor), which is called the \emph{second fundamental form} of $H$. The second fundamental form of the submanifold $H$ is denoted by $B$ and for every $x\in H$ is a map $B:T_{x}H\times T_{x}H\to T_{x}H^{\perp}$ defined as:
\begin{eqnarray*}
B(X,Y)=\nabla_{X}Y-(\nabla_{X}Y)^{T},
\end{eqnarray*}
where $(\nabla_{X}Y)^{T}$ is the orthogonal projection of the vector $\nabla_{X}Y$ into $T_{x}H$. It is easy to show that $B(X,Y)$ depends only on the value of $Y(x)$ and not on a choice of the field $Y$. The \emph{mean curvature} vector of the submanifold $H$ at the point $x\in H$ is denoted by $M_x$ and is defined by the relation
\begin{eqnarray*}
M=tr(B)=\frac{1}{k}\sum_{i=1}^{k}B(e_i,e_i),
\end{eqnarray*}
where $\{e_1,\cdots,e_k\}$ is an arbitrary orthonormal basis of $T_{x}H$. The trace is calculated with respect to the first fundamental form (metric) of $M$ restricted to $H$ (indeed since in the definition the choice of orthonormal basis of $T_{x}H$ depends on the metric of $M$). A submanifold  $H$ is called minimal if the mean curvature vector vanishes everywhere on $H$. 

We need to generalise the definition of conjugate points and cut-locus for the normal exponential maps with respect to submanifolds.
\begin{de}[Focal points and cut-focal Points]
A point $x\in M$ is called a \emph{focal} point of the submanifold $H$ if $exp_N$ is singular on $exp_N^{-1}(x)$.
A cut-focal point along a geodesic $\sigma$ orthogonal to the submanifold $H$ is a point on $\sigma$ such that the distance to the submanifold is no longer minimised along $\sigma$.
\end{de}
For every $(h,u)$ where $h\in H$ and $u\in H_h^{\perp}$ with $\Vert u\Vert=1$ define $e_c(h,u)$ to be the supremum over all $t>0$ such that $d(exp_N(h,tu),H)=t$.

Define $O_H\subset N_H$ to be the set of all $(h,tu)\in N_H$ with $\Vert u\Vert=1$ and $0\leq t<e_c(h,u)$.

Keeping this in mind, an immediate and usefull theorem is:
\begin{theo} \label{start}
Let $H$ be a $k$-dimensional compact submanifold of a compact $n$-dimensional Riemannian manifold $M$. Then we have
\begin{eqnarray*}
vol_n(M)&=&vol_n(O_H) \\
        &=&\int_{O_H}exp_N^{*}(dv_{M}) \\
        &=&\int_{H}\int_{\mathbb{S}^{n-k-1}}\int_{0}^{e_c(h,u)}t^{n-k-1}J_u(t)dt\,du\,dH,
\end{eqnarray*}
where $dv_{M}$ is the volume element in $M$.
\end{theo}

Let $M^n$ be a closed Riemannian manifold such that $K(M^n)\geq \delta>0$ (remember $K$ denotes the sectional curvature). Let $H^{k}\subset M^n$ be a closed $k$-dimensional submanifold such that
\begin{eqnarray*}
\Vert M\Vert(H^{k})\leq \kappa,
\end{eqnarray*}
where $\Vert M\Vert$ denotes the norm (calculated with respect to the metric of $M^n$) of the mean curvature vectorfield. 

The aim of this section is to prove the following:
\begin{theo} \label{sho}
Let $M^n$ and $H^{k}$ be as described above, then:
\begin{eqnarray*}
\frac{vol_{k}(H^{k})}{vol_n(M^n)}\geq \frac{vol_{k}(\mathbb{S}^{k}(\delta+\frac{\kappa^2}{k^2}))}{vol_n(\mathbb{S}^n(\delta))},
\end{eqnarray*}
where $\mathbb{S}^n(\delta)$ is the (model) space of constant sectional curvature equal to $\delta$.
\end{theo}

\emph{Remark}: Theorem \ref{sho} was proven by Heintze and Karcher in \cite{hk}. The reader can omit the proof and move to the next section. However, I shall present a proof for this Theorem which is slightly different from the proof in \cite{hk}.

I shall postpone the proof until the end of this section, and begin here to prove two key results for the proof of Theorem \ref{sho}. 

\begin{lem} \label{ma}
For any compact $k$-dimensional submanifold $H$ of a compact Riemannian manifold $M$ we have
\begin{eqnarray*}
\int_{\mathbb{S}^{n-k-1}}\int_{0}^{x(u)}(\frac{\sin(t\sqrt{\delta})}{\sqrt{\delta}})^{n-k-1}(\cos(t\delta)-\frac{\sin(t\delta)}{k\sqrt{\delta}}g(M,u))^k dt\,du \\
            = \frac{vol_n(\mathbb{S}^n(\delta))}{vol_k(\mathbb{S}^k(\delta+\frac{(\Vert M\Vert)^2}{k^2}))},
\end{eqnarray*}
where $x(u)$ is the first zero of $J_u(t)$, the Jacobian of the normal exponential map along $t\to exp_n(h,tu)$ in the space of constant curvature i.e. $\cos(x(u)\sqrt{\delta})=\frac{\sin(x(u)\sqrt{\delta})}{\sqrt{\delta}}g(M,u)$ and where for every $x\in H$, $u\in\mathbb{S}^{n-k-1}$ is a unit vector belonging to the unit sphere of the vector space $T_{x}H^{\perp}$ orthogonal of $T_{x}H$.
\end{lem}

\emph{Proof of Lemma \ref{ma}}

First let us look at the spheres of constant curvature. If $H$ is a closed submanifold of $\mathbb{S}^n(\delta)$, for sufficiently small $r$, the Jacobian of the map $exp_N$ is given by:
\begin{eqnarray*}
J_u(t)=(\frac{\sin(t\sqrt{\delta})}{t\sqrt{\delta}})^{n-k-1}.det(\cos(t\sqrt{\delta})I-\frac{\sin(t\sqrt{\delta})}{t\sqrt{\delta}}T_u),
\end{eqnarray*}
where $T_u$ is the Weingarten map of the second fundamental form of $H$ i.e. for $h\in H$ and $u\in TH_h^{\perp}$ $T_u:TH_h\to TH_h$ we have $g(T_u(v),w)=B(v,w)u$. Let $T_u$ be the Weingarten map of the embedding of the sphere $\mathbb{S}^k(\delta+\frac{\kappa^2}{k^2})$ into the sphere $\mathbb{S}^n(\delta)$. Since a sphere of latitude is totally umbilic, we have 
\begin{eqnarray*}
T_u=\frac{1}{k}g(M,u)I.
\end{eqnarray*}
Hence the radius of the sphere of latitude is equal to
\begin{eqnarray*}
\frac{1}{\sqrt{\delta+\Vert H\Vert^2/k^2}}=\frac{1}{\sqrt{\delta+\kappa^2/k^2}}.
\end{eqnarray*}
This shows that $\Vert M\Vert^2=\kappa^2$. Then, the Jacobian of the normal exponential map is equal to
\begin{eqnarray*}
J_u(t)=(\frac{\sin(t\sqrt{\delta})}{\sqrt{\delta}})^{n-k-1}(\cos(t\delta)-\frac{\sin(t\delta)}{k\sqrt{\delta}}g(M,u))^k
\end{eqnarray*}
where $x(u)$ is the first zero of $J_u(t)$ along $t\to exp_{\nu}(p,tu)$. This shows that both $J_u(t)$ and $x(u)$ are independant of the point $p$. Writing down the integral of the volume with respect to the Fermi coordinates defined by the sphere of lattitude, we get the desired result for the spheres.

For general manifolds, it is sufficient to notice that the calculation of the integral for the spheres is a pointwise computation, which would agree for every manifold. It equals the special case of the spheres and thus finishes the proof of Lemma \ref{ma}.
\begin{flushright}
$\Box$
\end{flushright}

\begin{lem} \label{do}
Let $H$ be a compact $k$-dimensional submanifold of a compact $n$-dimensional Riemannian manifold $X$ whose sectional curvature satisfies $K\geq \delta>0$. Then
\begin{eqnarray*}
vol_n(X) \leq  \left(\int_{H}\left(\int_{\mathbb{S}^{n-k-1}}\left(\int_{0}^{x(u)}(\frac{\sin(t\sqrt{\delta})}{\sqrt{\delta}})^{n-k-1}(\cos(t\delta)-\frac{\sin(t\delta)}{k\sqrt{\delta}}g(H,u))^k dt\right)du\right)dH\right).
\end{eqnarray*}
\end{lem}

\emph{Proof of Lemma \ref{do}}

We have
\begin{eqnarray*}
vol_n(X) \leq  \left(\int_{H}\left(\int_{\mathbb{S}^{n-k-1}}\left(\int_{0}^{e_{c}(h,u)}(\frac{\sin(t\sqrt{\delta})}{\sqrt{\delta}})^{n-k-1}.det(\cos(t\sqrt{\delta})I-\frac{\sin(t\delta)}{\sqrt{\delta}}T_{u})dt\right)du\right)dH\right) \\
\leq \left(\int_{H}\left(\int_{\mathbb{S}^{n-k-1}}\left(\int_{0}^{e_c(h,u)}(\frac{\sin(t\sqrt{\delta})}{\sqrt{\delta}})^{n-k-1}(\cos(t\delta)-\frac{\sin(t\delta)}{k\sqrt{\delta}}g(H,u))^k dt\right)du\right)dH\right) .
\end{eqnarray*}
Note that the integral of the above expression remains nonegative if we replace $e_c(h,u)$ by $x(u)$. This ends the proof of Lemma \ref{do}.
\begin{flushright}
$\Box$
\end{flushright}
\subsection{Proof of Theorem \ref{sho}}
Since the norm of mean-curvature of $H$ is supposed to be at least equal to $\kappa$ and according to the Lemmas \ref{ma} and \ref{do} we get: 
\begin{eqnarray*}
vol_n(X)\leq vol_n(\mathbb{S}^n(\lambda))\int_{H} \frac{1}{vol_{n-k}(\mathbb{S}^{n-k}(\lambda+\frac{\Vert H\Vert^2}{(n-k)^2})}dH.
\end{eqnarray*}

Thus:

\begin{eqnarray*}
\frac{vol_{n-k}(H)}{vol_n(X)}\geq \frac{vol_{n-k}(\mathbb{S}^{n-k}(\lambda+\frac{\kappa^2}{(n-k)^2}))}{vol_n(\mathbb{S}^n(\lambda))}.
\end{eqnarray*}
This proves Theorem \ref{sho}.
\begin{flushright}
$\Box$
\end{flushright}

Next section concerns Almgren-Pitts Min-Max theory.

\section{Almgren-Pitts Min-Max Theory}

As it was mentioned in Section $2$ of this paper, the invariant waist can be defined via a Min-Max procedure. In order to find a lower bound for the waist, it is then natural to study a variational problem on the \emph{space} of family of cycles sweeping out the manifold under study. These kinds of variational problems have been studied in the branch of geometric measure theory since the early sixties. The main idea of this theory, which will be useful here, can be roughly summarised in the following sentence:

\emph{For every family of cycles sweeping out a given closed Riemannian manifold $M$, call the cycle(s) with largest volume (Hausdorff measure) among other cycles (of this family), a critical cycle. Suppose $X$ sweeped out by a family of $k$-cycles is denoted by $\Pi_1$. There exists a sweep out $\Pi_2$ of $X$ by a family of $k$-cycles such that $\Pi_1$ is homotopic to $\Pi_2$ and the volume of the critical cycle in $\Pi_1$ is larger than the volume of the critical cycle in $\Pi_2$. Furtheremore, the critical cycle of $\Pi_2$ is a generalised minimal submanifold.}

The above is one of the main results of what is called \emph{the Almgren-Pitts Min-Max theory}. This section is devoted to this theory, more specifically, to giving the above sentence rigourus mathematical sense (which will be Theorem \ref{minmax}, announced in this section).
  
Chronologically speaking, shortly after his PhD thesis, Almgren developed a Morse Theory on the space of cycles in order to show the existence of some \emph{generalised} minimal submanifolds (subvarieties) in every (co)-dimension on any closed Riemannian manifold. The classical Morse Theory studies the length (energy) functional(s) on the (infinite dimensional) loop space of a Riemannian manifold. The Almgren-Morse Theory studies the volume functional on the (infinite dimensional) space of cycles. This idea, indeed, went back to other mathematicians-notably Birkhoff. Birkhoff's goal was to prove the existence of simple closed geodesics on surfaces diffeomorphic to the $2$-sphere. In order to achieve this goal, he showed that for any $1$-parameter family of closed curves sweeping out this surface, the longest curve always has a length greater or equal to the length of a (non-trivial) closed geodesic. This was the beginning of Min-Max type arguments in (geometric) variational problems. Almgren's Morse Theory aimed to generalise this type of argument in every dimension and codimension. However, Almgren's theory was never published (although recently Larry Guth kindly sent Almgren's early publication to me, and if I am not wrong, there is a chance that this old manuscript will be published, thanks to Guth). It was only in 1981 that Jon T. Pitts published a book on this matter. In this book \cite{pitts}, Pitts develops (and extends) Almgren's Theory and presents many results concerning the existence and regularity of minimal subvarieties. One important (and quite general) result proved in \cite{pitts} which extends Birkhoff's result to every dimension and co-dimension is the following:
\begin{theo}[Pitts-Almgren] \label{pa}
For each $k\leq n$, every compact $n$-dimensional Riemannian manifold of class $4$ supports a nonzero $k$-dimensional \emph{stationary integral varifold} which at each point in the manifold, is almost minimizing in all small neighborhoods of that point.
Morever if $3\leq n\leq 6$ and if the manifold is sufficiently smooth then it supports a nonempty closed embedded minimal submanifold of codimension $1$.
\end{theo}

As we can observe, there are some technical terms used in Theorem \ref{pa} which require some clarification. Additionaly, this theorem is a (huge) generalisation of Birkhoff's Theorem on the existence of simple closed geodesics on  a $2$-dimensional sphere. The interesting point is that the proof of Theorem \ref{pa} also generalises the Min-Max arguments of Birkhoff's proof.

It is not hard to guess that unfortunately, the result of a Min-Max procedure does not always give a smooth object and we should expect to work with objects other than smooth submanifolds. In every branch of mathematics which deals with variational problems, one begins with a class of objects, say $M$ (for example a space of functions, a space of currents, a space of geometric cycles, and etc.) studies the variational problem on this class $M$, and realises that the critical point(s) of the problem do not belong in $M$. This is why the Sobolev spaces of Levi-Fubini-Sobolev and the distributions of Schwartz exist.

Apart from the regularity issue, other important issues can be encountered as well. For instance, what is the best suited space of cycles on which the variational problem has to be studied? What should the topology on such a space be, so that the volume (Hausdorff measure) functional is continuous?

I shall give some quick definitions from geometric measure theory. The major geometric measure theoretical objects will be the \emph{varifolds} of Almgren. They are objects best suited for variational studies on the space of cycles. They are a replacement of submanifolds and can contain singularities. I will recall the main properties of the varifolds, then move to a major topological property of the space of cycles which will be Theorem \ref{hom}. Then we will have all the necessary material to announce the main Theorem of Almgren-Pitts Min-Max Theory, giving a rigourus sense to the vulgarisations of the beginning of this section. This will be Theorem \ref{minmax}. Finally, I will come back to the main concern of this paper, which is a lower bound for the invariant waist, and connect Almgren-Pitts Min-Max Theory to the waist by proving the Lemma \ref{critic} (which will be crucial in order to prove the main theorem of this paper: Theorem \ref{main}).

\begin{de}[Rectifiable sets]
A set $E\subset \mathbb{R}^N$ is called $k$-rectifiable if $H^k(E)<\infty$ and $H^k$-almost all of $E$ is contained in thge union of the images of countably many Lipschitz maps from $\mathbb{R}^k\to\mathbb{R}^N$.
\end{de}
In \cite{mor}, it is shown that a measurable set $S$ is rectifiable if and only if $H^k(S)<\infty$ and $H^k$-almost all of $S$ is contained in a countable union of $C^1$ embedded manifolds.

\begin{de}[Varifolds in $\mathbb{R}^N$]
An $k$-dimensional varifold $V$ in $\mathbb{R}^N$ is a Radon measure $V$ over $\mathbb{R}^n\times G(N,k)$, where $G(N,k)$ is the Grassmanian space consisting of the $k$-dimensional vector subspaces of $\mathbb{R}^N$. The space of $k$-dimensional varifolds in $\mathbb{R}^N$ is denoted by $V_k(\mathbb{R}^N)$ and the topology on this space is given by the weak topology of convergence of Radon measures.
\end{de}

\begin{de}[Mass of a varifold and variation measure]
The mass of a varifold $V$ is defined by
\begin{eqnarray*}
M(V)=V(\mathbb{R}^N\times G(N,k)).
\end{eqnarray*}
A varifold $V\in V_k(\mathbb{R}^N)$ defines a Radon measure on $\mathbb{R}^n$ which is denoted by $\Vert V\Vert$ and for every open set $U\subset \mathbb{R}^N$ is defined by
\begin{eqnarray*}
\Vert V\Vert(U)=V(U\times G(N,k)).
\end{eqnarray*}
\end{de}

\begin{de}[Rectifiable and Integral Varifolds]
A rectifiable varifold $V$ is a varifold induced by a rectifiable set $E$ and a Borel function $f:E\to\mathbb{R}_{+}$ defined by
\begin{eqnarray*}
\int_{\mathbb{R}^N\times G(N,k)}\phi(x,\pi)dV(x,\pi)=\int_{E}f(x)\phi(x,T_xE)dH^k(x),
\end{eqnarray*}
for every $\phi\in C_c(\mathbb{R}^N\times G(N,k))$, where $T_xE$ is the tangent cone of $E$ at the point $x$.
If $f$ is integer valued, we say that the varifold $V$ is integral Varifold.
\end{de}

\begin{de}[Mapping Varifolds]
Let $f:\mathbb{R}^N\to\mathbb{R}^N$ be a continuously differentiable map. Let $V\in V_k(\mathbb{R}^N)$, then $f_{\sharp}(V)\in V_k(\mathbb{R}^N)$ is the varifold defined by
\begin{eqnarray*}
\int_{\mathbb{R}^N\times G(N,k)}\phi(y,\sigma)d(f_{\sharp}V(y,\sigma))=\int_{\mathbb{R}^N\times G(N,k)}Jf(x,\pi)\phi(f(x),df_x(\pi))dV(x,\pi),
\end{eqnarray*}
where $Jf(x,\pi)$ is the Jacobian (determinant) of the differential $df_x$ restricted to the plane $\pi$.
\end{de}

\begin{de}[First Variation of Varifolds, Stationary Varifolds]
Let $V\in V_k(\mathbb{R}^N)$ and let $\chi(\mathbb{R}^N)$ be the space of smooth vector fields in $\mathbb{R}^N$ . The first variation of $V$ is the vector valued distribution $\delta V:\chi(\mathbb{R}^N)\to\mathbb{R}$ defined by:
\begin{eqnarray*}
\delta v(g)=\frac{d}{dt}(\Vert \psi(t,.)_{\sharp}V\Vert))\mid_{t=0}
\end{eqnarray*}
which is the initial rate of change of mass associated with the isotopy $\frac{\partial \psi}{\partial t}=g(\psi)$.

A varifold $V$ is said to be \emph{stationary}, if $\delta V(g)=0$ for every $g\in \chi(\mathbb{R}^N)$.
\end{de}

One should keep in mind that stationary varifolds play the same role as the minimal submanifolds.

So far, I have only defined varifolds in Euclidean spaces. How can we define these objects on Riemannian manifolds?

Varifolds on Riemannian manifolds are defined through Euclidean spaces. Therefore, we isometrically embed our Riemannian manifold in a certain Euclidean space (no need to remind the reader that this is always possible via Nash's famous Isometric Embedding Theorem). In such a space, we can define varifolds, and we simply ask that the support of such a varifold be contained in an (relatively) open subset of it. 

In Section $2$ of this paper, I defined the space of cycles and also presented a topological property to be satisfied in order for a family of cycles \emph{sweeping out} a given metric space. The topology of the space of cycles plays an important role in the geometry of the Min-Max. The topology we put on the space of cycles was the $\flat$-norm. However, it is important to point out that there are other topologies one could define on this space. An important topology which makes the volume functional a continuous functional (and is considered in \cite{pitts}) is called the $M$-topology. In \cite{alm}, Almgren proves a generalised Dold-Thom Theorem, computing the homotopy groups of the space of cycles and relating them to the homology of the manifold on which the cycles are supported. In \cite{pitts}, Pitts (slightly) extends Almgren's proof, and together with ideas in \cite{alm} proves the following important theorem: 
\begin{theo}[Almgren-Pitts] \label{hom}
Let $X$ be a closed $n$-dimensional Riemannian manifold and let $1\leq k\leq n$. Let $G$ be the coefficient group of a classical (singular) homology theory $H$. Then
\begin{equation}
H_{m+k}(X,G)\cong \pi_{m}(C_k(X,G),\{0\}), \label{eq:solve}
\end{equation}
where $\pi_{m}(C_k(X,G),\{0\})$ is the $m$th homotopy group of the space of $k$-cycles of the manifold $X$ with basepoint the trivial $k$-cycle $0$. The isomorphism in (\ref{eq:solve}) is natural. This result remains valid if one replace the $\flat$ topology of the space of cycles by the $M$-topology.
\end{theo}

According to Theorem \ref{hom}, every homotopy class of mappings in the space $C_k(X,G)$ naturally defines a min-max procedure. The min-max procedure related to the sweep out of the manifold $X$ corresponds to the $(n-k)$th-homotopy group of the space $C_k(X,G)$. The variational problem to any min-max procedure can simply be defined as follows:
\begin{de}[Variational Problem Related to Min-Max]
For every homotopy class of 

mappings $[\Pi]$ in $C_k(X,G)$, define:
\begin{eqnarray*}
L([\Pi])=\inf_{\Pi'\simeq \Pi}(\sup_{\pi\in \Pi}vol(\pi)),
\end{eqnarray*}
where $vol$ denotes the Hausdorff measure and $\simeq$ denotes the homotopy relation.
\end{de}





We are ready to state the main theorem of Almgren-Pitts min-max theory. The proof is very long and technical and can be found in \cite{pitts}:
\begin{theo} \label{minmax}
Let $X$ be a compact Riemannian manifold of dimension $n$. Suppose $X$ be isometrically embedded in $\mathbb{R}^N$. let $1\leq k\leq n$. If $[\Pi]$ is an homotopy class of mappings into $C_k(X,G)$, then there exists a stationary integral varifold $V\in V_k(X)$ such that:
\begin{eqnarray*}
\Vert V\Vert(\mathbb{R}^N)=L(\Pi).
\end{eqnarray*}
Furtheremore, for each $x\in X$, there exists a positive number $r$ such that $V$ is $G$ almost minimisng in $X\cap A_N(x,s,r)$ for all $0<s<r$, where $A_N(x,s,r)$ denotes the annulus in $\mathbb{R}^N$ centered at $x$ with the corresponding radii. And for the varifold $V$ to be $G$ almost minimising means that it can be approximated by $G$-mass minimising currents.
\end{theo}

With Theorem \ref{minmax} to hand, we can come back to our main interest of giving a lower bound for the waists of a closed Riemannian manifold. At this stage, it is natural to ask the following questions:
\begin{itemize}
\item \emph{Does a lower bound on the mass of stationary integral $k$-dimensional varifolds give a lower bound on the $k$-waist}?
\item \emph{If so, how sharp is the result (i.e. a sharp lower bound on the mass of minimal subvarieties yields to a sharp lower bound on the waist})?
\end{itemize}

the next lemma answers the first question and relates Almgren-Pitts Min-Max Theory to the waist invariant:

\begin{lem} \label{critic}
Let $X$ be a closed Riemannian $n$-manifold. Let $1\leq k\leq n$ and suppose a $w_0\in \mathbb{R}_{+}$ exists such that for every $k$-dimensional stationary integral varifold (or generalised minimal submanifold) $H_k$, we have
\begin{eqnarray*}
M(H_k)\geq w_0.
\end{eqnarray*}
Then for every smooth map $f:X\to\mathbb{R}^{n-k}$, there exists a point $z\in\mathbb{R}^{n-k}$ such that
\begin{eqnarray*}
vol_k(f^{-1}(z))\geq w_0.
\end{eqnarray*}
In other words, we have
\begin{eqnarray*}
wst_k(X)\geq w_0.
\end{eqnarray*}
\end{lem}
\emph{Proof of Lemma \ref{critic}}:

Suppose the manifold $X$ is isometrically embedded in $\mathbb{R}^N$. Let $f:X\to \mathbb{R}^{n-k}$ be a smooth map. The map $f$ defines a continuous map $F:\mathbb{R}^{n-k}\to C_k(X,G)$ (for a coefficient group which can be choosen to be $\mathbb{Z}_2$). By a slight perturbation (if necessary), we can assume that there exists a point $w\in\mathbb{R}^{n-k}$ such that $f^{-1}(w)$ is a trivial cycle denoted by $*$. The unit cube $I^{n-k}$ being homeomorphic to $\mathbb{R}^{n-k}$, we chose a homeomorphism $h$ of the pair $(I^{n-k},\partial I^{n-k})$ to $(\mathbb{R}^{n-k},w)$, where $\partial I^{n-k}$ is the boundary of the closure of the unit cube $I^{n-k}$. Hence 

\[
\begin{tikzcd}
(I^{n-k},\partial I^{n-k}) \arrow{r}{\pi_{E}} \arrow[swap]{dr}{G} &(\mathbb{R}^{n-k},w)   \arrow{d}{\tilde{F}} \\
& (C_k(X,G),*) 
\end{tikzcd}
\]

we get a continuous map $G:(I^{n-k},\partial I^{n-k})\to (C_k(X,G),*)$. According to Theorem \ref{hom}, the homotopy class of the mappings $[G]$, defines a min-max procedure and according to Min-Max Theorem \ref{minmax}, there exists a family of $k$-cycles, homotopic to the family $G$ and a stationary integral (almost minimising) $k$-varifold $V$ in $X$ such that:
\begin{eqnarray*}
\sup_{z\in\mathbb{R}^{n-k}}vol_k(f^{-1}(z))&\geq& \Vert V\Vert(\mathbb{R}^N) \\
                                           &=& M(V)\\
                                           &\geq& w_0.
\end{eqnarray*}

This ends the proof of Lemma \ref{critic}.
\begin{flushright}
$\Box$
\end{flushright}

The method to keep in mind for one who is interested in finding a (lower) bound on the $k$-waist of a closed Riemannian manifold, is to try to find a lower bound for the mass (volume) of all stationary integral $k$-varifolds which are combinatorialy a $k$-cycle. In the last section, I will discuss the sharpness of the result of Lemma \ref{critic}. We will see that even if this lemma has provided us with a satisfactory algorithm to find a lower bound for the waist of a closed Riemannian manifold, simple examples exist showing that this result can be far from sharp.

\section{Waist of Positively Curved Riemannian Manifolds and the Proof of Theorem \ref{main}}

We are now set to give a complete proof of Theorem \ref{main}. We assume $X$ to be a closed Riemannian manifold of dimension $n$ such that  $K(X)\geq \kappa>0$. According to Lemma \ref{critic}, in order to give a lower bound for the $k$-waist of $X$, it is sufficient to give a lower bound for the mass of stationary integral $k$-varifolds of $X$. Let $Z_k$ be a stationary integral varifold in $X$. If $Z_k$ is a minimal submanifold, we apply directly the result of Theorem \ref{sho} and get:
\begin{eqnarray*}
\frac{vol_k(Sup_{z\in\mathbb{R}^k}vol_{n-k}(f^{-1}(z)))}{vol_n(X)}&\geq& \frac{vol_k(Z_k)}{vol_n(X)}\\
                                                               &\geq&\frac{vol_{k}(\mathbb{S}^{k}(\kappa))}{vol_n(\mathbb{S}^n(\kappa))}.
\end{eqnarray*}

And thus the required result is proven. Of course, as it was mentioned before, $Z_k$ can have complicated singularities, and for this case, we can no longer apply directly the result of Theorem \ref{sho}. The complexity of the geometry and structure of the singularities of such an object remains unclear. Through deep regularity results of Allard in \cite{alla}, together with higher differentiability theory of Morrey (see \cite{mor} and \cite{pitts}), we only can conclude that an open dense subset of the support of this varifold corresponds to a smooth $k$-dimensional submanifold with everywhere mean-curvature equal to zero. Unfortunately, this is still not good enough for us to give a lower bound for the size (Hausdorff measure) of this varifold. One could try a smoothing procedure and approximate this varifold by smooth submanifolds, but in high codimension there is not a good control on the mean-curvature of the submnaifolds approximating this stationnary integral varifold (this was pointed out to me by Camillo De Lellis). In order to solve the issue of the singularities and provide a (sharp) lower bound for the size of a $k$-dimensional almost-minimising stationary integral varifold, I will use a standard procedure. When we deal with stationary integral varifolds with possible singularities, this procedure is called the \emph{monotonicity formula} which will be explained here:

I begin first by proving a localised version of Theorem \ref{sho}:

\begin{lem} \label{local}
Let $H$ be a stationnary integral varifold of dimension $k$ in $X$. Suppose $K(X)\geq \delta$. Let $x$ be any regular point in $H$. Let $r>0$. Then
\begin{eqnarray*}
\frac{vol_k(H\cap B(x,r))}{vol_n(B(x,r))}\geq \frac{vol_k(S^k(\kappa)}{vol_n(\mathbb{S}^n(\kappa)},
\end{eqnarray*}
where $B(x,r)$ is an intrinsic ball of radius $r$ centered at the point $x\in X$.
\end{lem}

\emph{Proof of Lemma \ref{local}}:

According to Lemma \ref{do} we have:
\begin{eqnarray*}
vol_n(B(x,r))&\leq& U_r(H\cap B(x,r)) \\
             &=&\left(\int_{H\cap B(x,r)}\left(\int_{\mathbb{S}^{n-k-1}}\left(\int_{0}^{r}(\frac{\sin(t\sqrt{\delta})}{\sqrt{\delta}})^{n-k-1}(\cos(t\delta) dt\right)du\right)dH\right)\\
             &\leq& \left(\int_{H\cap B(x,r)}\left(\int_{\mathbb{S}^{n-k-1}}\left(\int_{0}^{x(u)}(\frac{\sin(t\sqrt{\delta})}{\sqrt{\delta}})^{n-k-1}(\cos(t\delta) dt\right)du\right)dH\right).
\end{eqnarray*}
Now use Lemma \ref{ma} to prove the following:
\begin{eqnarray*}
\left(\int_{H\cap B(x,r)}\left(\int_{\mathbb{S}^{n-k-1}}\left(\int_{0}^{x(u)}(\frac{\sin(t\sqrt{\delta})}{\sqrt{\delta}})^{n-k-1}(\cos(t\delta) dt\right)du\right)dH\right)\leq \frac{vol_n(\mathbb{S}^n(\kappa))}{vol_k(\mathbb{S}^k(\kappa))}.vol_k(H\cap B(x,r)).
\end{eqnarray*}
Combining the above two equations, we get:
\begin{eqnarray*}
\frac{vol_k(H\cap B(x,r))}{vol_n(B(x,r))}\geq \frac{vol_k(S^k(\kappa)}{vol_n(\mathbb{S}^n(\kappa))},
\end{eqnarray*}
and this ends the proof of Lemma \ref{local}.
\begin{flushright}
$\Box$
\end{flushright}

\emph{Remark}:

The result of Lemma \ref{local} is far from being sharp. I believe a better lower bound could have been proven. Let $\bar B_n(.,r)$ denotes any ball of radius $r$ in $\mathbb{S}^n(\kappa)$. I believe that:
\begin{eqnarray*}
\frac{vol_k(H\cap B(x,r))}{vol_n(B(x,r))}\geq \frac{vol_k(S^k(\kappa)\cap \bar B_n(y,r))}{vol_n(\bar B_n(y,r))},
\end{eqnarray*}
where $y\in \mathbb{S}^k$.

Coming back to the proof of Theorem \ref{main}, since an open and dense subset of the support of $Z_k$ is smooth and has mean-curvature everywhere vanishing then applying  Lemma \ref{local} and for $r\to \infty$, the proof of Theorem \ref{main} immediately  follows.
\begin{flushright}
$\Box$
\end{flushright}

\section{Ideas Related to Conjecture \ref{gro}}

This section concerns Conjecture \ref{gro} for $\varepsilon>0$. At first glance, one might think it would be easy to prove Conjecture \ref{gro} directly from Theorem \ref{main} by simply integrating. Indeed it is, since the map $f$ is supposed to be smooth, the fiber $f^{-1}(z)$, maximising the volume among other fibers can be assumed to be smooth as well. Then, according to the tube formula, we have: 
\begin{equation}
vol_n(f^{-1}(z)+\varepsilon)=\int_{0}^{\varepsilon}vol_{n-1}(\partial U_r(f^{-1}(z)))dr. \label{eq:solvee}
\end{equation}
The aim would be to give a (sharp) lower bound for $vol_n(f^{-1}(z)+\varepsilon)$ which would result in a lower bound for the right hand side of the equation \ref{eq:solvee}. This could be possible if for every $0\leq t\leq \varepsilon$, we could find an upper bound for the supremum of the norm of the mean-curvature of the boundary of the (closed) tube $\partial U_t(f^{-1}(z))$. However this actually seems to be impossible even through the Riccati equations along geodesics normal to $f^{-1}(z)$.

This is the first issue. The second (and more fundamental) issue is the fact that nothing states that if $f^{-1}(z)$ maximises the volume over other fibers, it also maximises the other functional $vol_n(. +\varepsilon)$ over \emph{all} the fibers. Indeed for a given $0<\varepsilon$, the functional $vol_n(.+\varepsilon)$ is not continuous on the space of $k$-cycles $C_k(X,G)$ (with $\flat$-norm topology or for $M$-topology). Hence, one can not repeat the arguments of Almgren-Pitts Min-Max theory given in \cite{pitts} in order to prove a similar theorem such as Theorem \ref{minmax}. 

Similar to Section $4$, I define:
\begin{de}[Variational Problem Related to $\varepsilon$-Min-Max]
For every $\varepsilon>0$ and for every homotopy class of mappings $[\Pi]$ in $C_k(X,G)$, define:
\begin{eqnarray*}
L_{\varepsilon}([\Pi])=\inf_{\Pi'\simeq \Pi}(\sup_{\pi\in \Pi}vol_n((\pi)+\varepsilon)),
\end{eqnarray*}
where $vol_n$ denotes the Riemannian volume and $\simeq$ denotes the homotopy relation.
\end{de}
  
I believe the following Conjecture (at least for $\varepsilon>0$ small enough) to be true:





\begin{conj} \label{minmaxep}
Let $X$ be a compact Riemannian manifold of dimension $n$. Suppose $X$ is isometrically embedded in $\mathbb{R}^N$. Let $1\leq k\leq n$. If $[\Pi]$ is an homotopy class of mappings into $C_k(X,G)$, then a stationary integral varifold $V\in V_k(X)$ exists such that:
\begin{eqnarray*}
\Vert V\Vert(\mathbb{R}^N)=L_{\varepsilon}(\Pi).
\end{eqnarray*}
Furtheremore, for each $x\in X$, a positive number $r$ exists such that $V$ is $G$ almost minimisng in $X\cap A_N(x,s,r)$ for all $0<s<r$, where $A_N(x,s,r)$ denotes the annulus in $\mathbb{R}^N$ centered at $x$ with the corresponding radii.
\end{conj}

\emph{Questions}:
\begin{itemize}
\item How big a role does the value of $\varepsilon$ play in the truth of Conjecture \ref{minmaxep}?

\item If Conjecture \ref{minmaxep} is true, could we use the following deep theorem proven by Mahmoudi-Mazzeo-Packard in \cite{mazzeo}:
\begin{theo}[Mahmoudi-Mazzeo-Pacard] \label{maz}
Let $H^k\subset X^n$ be a closed (embedded or immersed) minimal submanifold (with possible singularities) which is nondegenrate in the sense that its Jacobi operator is invertible. Then for sufficiently small $\varepsilon>0$ and for every $0<r\leq\varepsilon$, $\partial U_r(H)$ may be perturbed to a constant mean curvature hypersurface $H_r^{n-1}$ with $\Vert M\Vert=\frac{n-k-1}{n-1}r^{-1}$. Furtheremore, the index of the hypersurface $H_r^{n-1}$ tends to $\infty$ as $r\to 0$.
\end{theo}

in order to prove a partial answer to Conjecture \ref{gro} for manifolds supporting $k$-dimensional minimal \emph{submanifolds}?

\item Is Theorem \ref{maz} true if one replaces minimal submnaifolds with stationnary integral varifolds?

\end{itemize}

\section{Sharp Waists and Minimal Submanifolds}

Note that according to Theorem \ref{sho}, it makes sense to believe that the optimal waist is supported on a minimal submanifold (i.e. a submanifold for which the mean-curvature vector vanishes everywhere). 

Of course the result of Theorem \ref{main} shows that if a positively-curved Riemannian manifold supports a codimension $k$ minimal submanifold $H$, then it's $k$-waist is bounded from below by the $vol_{n-k}(H)$ (for example in the case of constant curvature spheres). Let us take, as an example, the complex projective space $\mathbb{C}P^n$ induced with its \emph{standard} metric (the Fubini-Study metric). Let us look for the $k$-waist where $k$ is an odd number. Then, a good candidate for the waist would be a totally geodesic submanifold. Remark though that there is no odd-dimensional totally geodesic $H\subset \mathbb{C}P^n$ (of course where $dim H>1$). Since $\mathbb{C}P^n$ induced with its canonical Riemannian metric has $k\geq \frac{1}{4}$, then we could use our main Theorem \ref{main} and have a lower bound for its waist in every (co)dimension for $\varepsilon=0$.

\emph{Questions}: 
\begin{itemize}
\item How sharp is this result (for $\varepsilon=0$)?
\item Could we find a sharp lower bound for the waists of the Canonical $\mathbb{C}P^n$ without having to prove Conjecture \ref{minmaxep}, and by using different techniques?
\item Could the partition-type proof of the waist of the canonical sphere used in \cite{grwst} and \cite{memwst} be applied for $\mathbb{C}P^n$ (and/or for more general symmetric spaces)?

Thanks to R.Karasev, A.Hubard and B.Matschke, I learned that the partition arguments used in \cite{grwst} and \cite{memwst} can be generalised extensively to \emph{any} manifolds-provided the \emph{space} of the \emph{cutting elements} are large enough. For symmetric spaces such as $\mathbb{R}P^n$ and $\mathbb{C}P^n$, using arguments like in \cite{karasev} and \cite{matsh}, I believe for any continuous map $f:\mathbb{C}P^n\to\mathbb{R}^k$, where $k\leq n$, we can find a convex partition \emph{adapted} to this map. (for definition of adapted maps, consult \cite{memwst}). From here, we need to define a new class of convexely-derived measures for convex sets of $\mathbb{C}P^n$ in the same way as we did for the case of the round-sphere. A good candidate for this class could be the class of $x^k\sin^m$-concave measures for appropriate $k$ and $m$.
\item What about the waist of symmetric Riemmanian manifolds other than the real and complex projective spaces?
\end{itemize}

Using the techniques of this paper for finding lower bound on the waist of closed Riemannian manifolds with $\vert K\vert\leq 1$ would not lead to sharp results. The next picture shows a $2$-dimensional Riemannian manifold in which there is a big gap between the length of the \emph{waist} and the \emph{smallest} minimal submanifold:
\begin{center}
\includegraphics[width=5in]{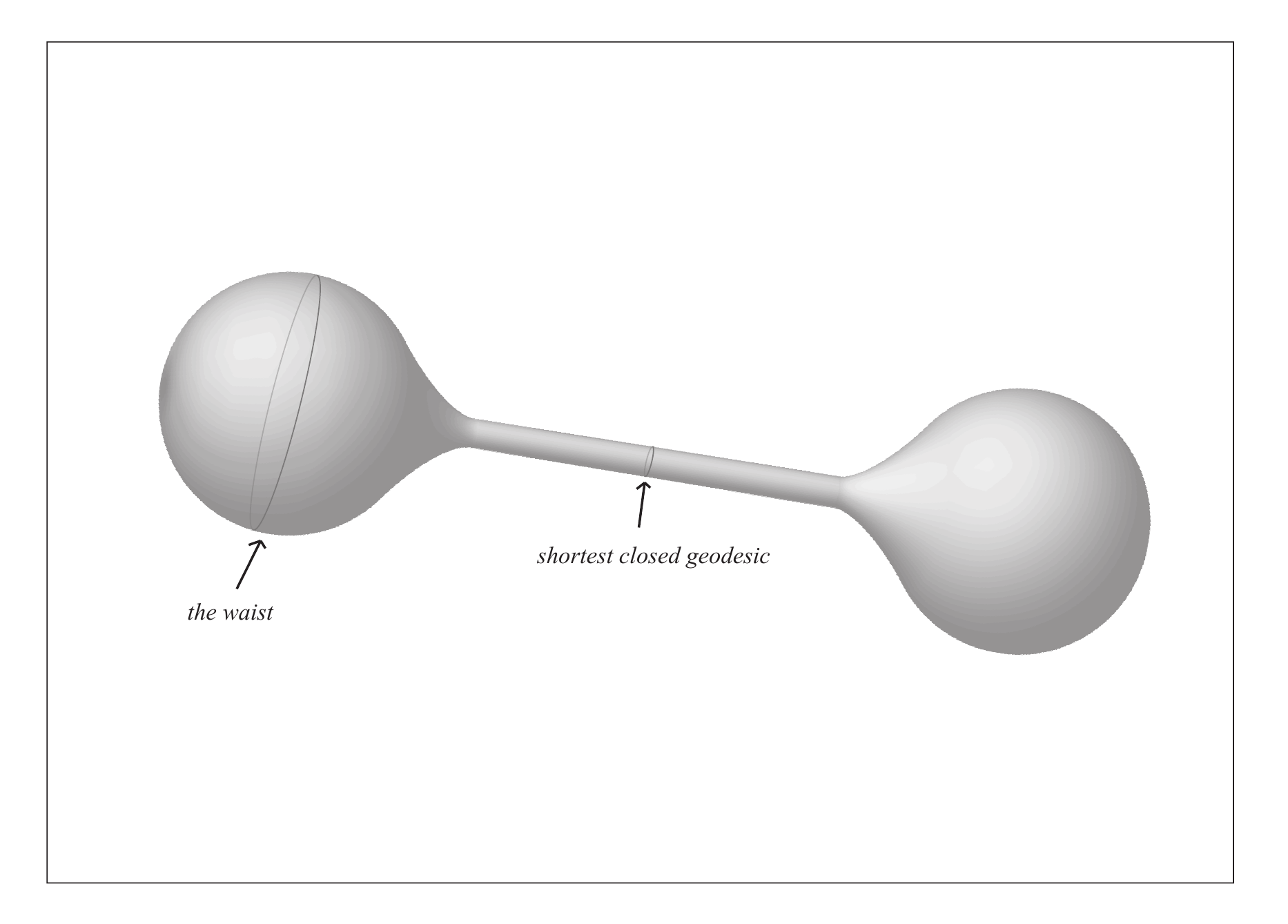}

Waist and Minimal Subvarieties
\end{center}

An ambitious and interesting research plan could be studying techniques in order to determine sharp waists of different metric-measure spaces. Even for the case in which the manifold has a large symmetry group, the sharp waist remains an open question.

\bibliographystyle{plain}
\bibliography{poscur}

\end{document}